\documentclass{amsart}
\usepackage{amssymb, amscd, amsmath, amsthm, epsf, epsfig, latexsym}

\newtheorem{theorem}{Theorem}
\newtheorem{corollary}{Corollary}[theorem]
\def\zz{{\bf Z}}
\def\calc{\mathcal{C}} 
\def\calg{\mathcal{G}}
\def\cala{\mathcal{A}} 
\def\zz{{\bf Z}}

\begin{document}
 
\title{Knot Signature Functions are Independent}

\author{Jae Choon Cha}
\author{Charles Livingston}

\address{Department of Mathematics, Indiana University, Bloomington, IN 47405}

\email{jccha@indiana.edu}
\email{livingst@indiana.edu}

 
\date{\today}

\begin{abstract}  A Seifert matrix is a square integral matrix $V$ satisfying
$\det(V - V^T) = \pm 1$.  To such a matrix and unit complex number $\omega$ there
is a signature, $\sigma_\omega(V) =
\mbox{sign}( (1 - \omega)V + (1 - \bar{\omega})V^T)$. Let $S$ denote the set of
unit complex numbers with positive imaginary part.  We show  
$\{\sigma_\omega\}_ {    \omega \in  S     }$   is linearly independent, viewed
as a set of  functions on the set of all Seifert matrices.

If $V$ is metabolic, then $\sigma_\omega(V) = 0$ unless $\omega$ is a root of the
Alexander polynomial, $\Delta_V(t) = \det(V - tV^T)$.  Let  $A$ denote the set of
all unit roots of all Alexander polynomials with positive imaginary part.  We
show that       $\{\sigma_\omega\}_ {   
\omega \in  A     }$ is linearly independent when viewed as a set of functions on
the set of all metabolic Seifert matrices.

To each knot $K \subset S^3$ one can associate a Seifert matrix $V_K$, and
$\sigma_\omega(V_K)$ induces a knot invariant.  Topological applications of our
results include a proof that the set of functions 
$\{\sigma_\omega\}_ {    \omega \in  S     }$ is linearly independent on the set
of all knots and that the set of two--sided averaged signature functions,
$\{\sigma^*_\omega\}_ {    \omega \in  S     }$, forms a linearly independent set
of homomorphisms on the knot concordance group.  Also, if $\nu$ is the root of
some Alexander polynomial, there is a slice knot
$K$ satisfying
$\sigma_\omega(K)
\ne 0$ if and only if $\omega = \nu$ or $\omega =\bar{\nu}$.
We demonstrate that the results extend to the higher
dimensional setting. 
\end{abstract}

\maketitle

 \section{Introduction}

Associated to each knot $K \subset S^3$ there is a {\it Seifert matrix} $V_K$. 
The set of such Seifert matrices consists of those square integral matrices $V$
satisfying $\det(V - V^T) = \pm 1$, where
$V^T$ denotes the transpose.  For each unit complex number  $\omega$  the {\it
hermitianized Seifert form},
$V_\omega$, is defined by 
$V_\omega  = (1 - \omega)V + (1- \bar{\omega})V^T$; the signature of this matrix 
 is denoted $\sigma_\omega(V)$.  It is possible to associate different Seifert
matrices to a given knot, however the value of $\sigma_\omega(V_K)$ is known to
depend only on
$K$ and not the choice of Seifert matrix \cite{le1, tr}; $\sigma_\omega(V_K)$ is
usually denoted
$\sigma_\omega(K)$. 

 Let $S$ denote the set of all unit complex numbers with positive imaginary
parts. In this paper we study the linear independence of the set of these
signature functions, $\{\sigma_\omega\}_\omega \in S$, viewed as real valued
functions on the set of all Seifert matrices, and consequently as functions on
the set of all knots.  Our first result is the following:

\begin{theorem}\label{thmind} The  set of functions, $\{\sigma_\omega\}_ {   
\omega \in  S     }$,  is linearly independent. 
\end{theorem}
 Previously the only sets $D$ for which it was known that  $\{\sigma_\omega\}_
{    \omega \in  D     }$ is linear independent were certain discrete subsets of
$S$, \cite {le1, tr}.  Applications of this result include the demonstration that a
number of results that hold in high-dimensional knot theory fail in dimension 3. 
This is briefly summarized in Section~\ref{appsect}.

A Seifert matrix is necessarily  of  even dimension, say $2g$. It is called {\it
metabolic} if there is a summand of dimension $g$ of $\zz^{2g}$ on which the
associated bilinear form vanishes.  It is called {\it hyperbolic} if $\zz^{2g}$
is the direct sum of two such summands.  Over the rational numbers these are
identical concepts, but that is not the case over $\zz$.  In particular,
$\sigma_\omega(V)$ is identically 0 if $V$ is hyperbolic, but can be nonzero for
a metabolic form $V$ if $\omega$ is a root of the {\it Alexander polynomial} of
$V$, $\Delta_V(t) = \det(V - tV^t)$. A detailed analysis of the signature
functions of metabolic forms was accomplished in \cite{le3}.  The topological
significance of these concepts is that metabolic forms correspond to {\it slice
knots} and hyperbolic forms correspond to {\it double null--concordant knots}
\cite{su}.  Renewed interest in double null-concordance, as summarized in
Section~\ref{appsect}, motivates our study of the signature functions associated to
metabolic Seifert forms.

Polynomials that occur as Alexander polynomials are precisely those polynomials
that are symmetric,
$\Delta(t^{-1}) = \pm t^j \Delta(t)$, for some $j$,  and satisfy $\Delta(-1) =
\pm 1$.  Let $D \subset S$ denote the set of unit roots of Alexander polynomials
with positive imaginary parts.

\begin{theorem}\label{thmmet} The   set of functions, $\{\sigma_\omega\}_ {   
\omega \in  D     }$,  is linearly independent on the set of all metabolic
Seifert matrices. 
\end{theorem}

This has the following topological corollary.

\begin{corollary}If $\nu$ is a unit root of some Alexander polynomial, there is a
slice knot $K$ with
  $\sigma_\omega(V_K) \ne 0$ if and only if $\omega = \nu \mbox{\ or\ }
\bar{\nu}$.

\end{corollary}

\section{Independence of Signature Functions}
\vskip.1in
\begin{proof}[Proof of Theorem \ref{thmind}] For a given Seifert matrix $V$,
$\sigma_\omega(V)$ can be viewed as an integer valued function of $\omega
\in S$.  Simple arguments show that jumps of this function can occur only at
those values of $\omega$ that are roots of   $\Delta_V(t)$, and if the root is
simple the jump  is nontrivial.  Also, for each
$V$, $\sigma_\omega(V)= 0$ for all
$\omega$ close to 1.  (For proofs, see, for instance,  \cite{le1}.)

We will next show that for any given $\omega \in S$ there is an $\omega' \in S$
arbitrarily close to
$\omega$ and a Seifert matrix $V$ whose signature function has its only
nontrivial jump at $\omega'$.  From this the theorem follows since one easily
constructs, for any finite set   $
\{\sigma_{\omega_i}\}_ {   1 \le i \le N}$ and any chosen $k, 1 \le k  \le N$, a
Seifert matrix
$V_k$ with
$\sigma_{\omega_k}(V_k)
\ne 0$ and
$\sigma_{\omega_i}(V_k) = 0$ if $1 \le i \le N$ and $  i \ne  k$. 

To construct the   desired matrix $V$ we construct a  polynomial $\Delta$ having
a unique
  root in $S$ at a point $\omega'$ that can be made as close to $\omega$ as
desired.  Since
$\Delta$ will be constructed to be integral, to satisfy $\Delta(1) = \pm 1$  and
to be symmetric, it is  the Alexander polynomial of some Seifert matrix $V$,
\cite{se}.

For a given $r$, $-1 < r <1$, consider the polynomial $$F_r(t) = (t-1)^2 (t^2 -
2rt +1) =  t^4+(-2-2 r) t^3+(4 r+2) t^2+(-2-2  r) t+1.$$ It is easily seen that
$F_r(t)$ has a pair of (unit) complex roots with real part $r$. Let $r =
$Re($\omega$) so that $F$ has roots at $\omega$ and $\bar{\omega}$. For a given
$\epsilon$ there is a
$\delta$ such that a perturbation of the coefficients of $F_r$ by less than
$\delta$ moves the roots less than
$\epsilon$.  Choose a rational approximation $a/b, $  $b > 0$, to $r$ so that
replacing $r$ in the coefficients of $F_r$ by $a/b$ changes the coefficients by
less than $\delta / 2$.  Furthermore, choose $b$ large enough so that $1 / b<
\delta / 2$.  Then the roots of
$$G(t) = t^4+(-2-2 {a \over b}) t^3+(4{a \over b}+2 - {1 \over b}) t^2+(-2-2 {a
\over b} )
 +1 $$ are within $\epsilon$ of those of $F_r$.  Multiplying through by $b$
yields the desired  polynomial:
$$ \Delta(t) = bt^4+(-2b-2 a) t^3+(4 a +2b - 1) t^2+(-2b-2 a) t+b .$$ 

Since 
$\Delta(1) = -1$ and $b > 0$, $\Delta(t)$ has at least two real roots.  The
remaining roots
$\omega'$ and $\bar{\omega}'$ are within $\epsilon$ of $\omega$ and
$\bar{\omega}$.   Since
$\Delta(t)$ is real, symmetric, and it has exactly two nonreal roots, these roots
must be lie on the unit circle. The result follows.  
 
\end{proof}
 \vskip.1in

\section{ A Remark  on Concordance}\label{concord} There is an equivalence relation
on the set of Seifert forms called {\it algebraic concordance}: $V_1$ and $V_2$ are
{\it concordant} if $V_1 \oplus -V_2$ is metabolic.  The  set of equivalence classes
forms a group,
$\calg$, the {\it algebraic concordance group}. (See \cite{le1}.)  The signature
functions are not well--defined on
$\calg$, but this difficulty can be overcome by considering    the averaged
signature function,
$\sigma^*_\omega$, obtained as the two--sided average of
$\sigma_\omega$.  Our theorem is easily seen to apply to the set
$\{\sigma^*_\omega\}_{\omega \in S}$, where this is now a set of {\it
homomorphisms}, not simply functions, on $\calg$.

The same result holds for the corresponding topological construction, the {\it
concordance group of knots}, $\calc$.  Further details will be summarized in the
final section.

\vskip.2in

\section{Signature Functions of Metabolic Forms}

A simple algebraic calculation shows that the matrix $  V_\omega$ is  singular
precisely when the unit
$\omega$ is a root of  
$\Delta_V(t)$.  If $V$ is metabolic  it follows readily that $\sigma_\omega(V) =
0$, except perhaps when $\omega$ is a root of $\Delta_V(t)$.  In \cite{le3},
Levine
  constructed metabolic matrices  with nontrivial signature at   $\omega = e^{ 
\pi i / 3}$, the root of $t^2 - t +1$.

Our goal is to show that Levine's construction can be expanded to cover any root
of any Alexander polynomial.  One slightly subtle point in the proof of Theorem
\ref{thmmet} is in dealing with Alexander polynomials that have several unit
roots; in such cases we must be able to specify at exactly which roots the
signature function is nonzero.

\begin{proof}[Proof of Theorem \ref{thmmet}] Let $\Delta(t) = \sum_{i = 0}^{2g}
d_i t^i$   be an Alexander polynomial.  Since multiplication by $\pm t^j$ does
not change the unit roots of
$\Delta(t)$, we can assume that
$d_0
\ne 0$,  $d_{2g -i} = d_i$ and that $\Delta(1) = 1$.  

Consider the matrix $V$ below, with $0_{2g}$ a $2g \times 2g$ matrix of zeroes
and $I_g$ the $g
\times g$ identity.  The $g \times g$ and $2g \times 2g$ integer matrices $A_g$
and $B_{2g}$ will be specified in the course of the proof.  We will choose
$B_{2g}$ to be symmetric, so assume so throughout the discussion.

$$ V =
\left(
\begin{tabular}{ c   c}
 {\huge{0}}$_{2g}$& 

$\left(  \begin{tabular}{c   c}                    I$_g$ & $A_g$ \\ 
 I$_g$ & 0$_g$ \\ \end{tabular}
\right) $

 \\  

$\left(  \begin{tabular}{c   c}                    0 & I$_g$ \\ 
 $A_{g} ^T$ & I$_g$ \\ \end{tabular}
\right) $

 &   \mbox{\LARGE $  \Omega    B_{2g}$}  \\
\end{tabular}
\right)
$$

Since $\det(V - V^T) = 1$, $V$ is the Seifert matrix. The half--dimensional block
of zeroes implies that $V$ is metabolic.  

Simple work with the matrices and algebraic manipulations yield  that
$\Delta_V(t)$ is the square of
$\det ((1 -t)^2 A_g +t)$, or, more usefully, 
$$\Delta_V(t) = \left( (1-t)^{2g} \lambda(\frac{t}{(1-t)^2} )\right)^2,$$ where
$\lambda(x) =\det(A_g + xI_g$).  Writing
$\lambda(x) = \sum_{j=0}^{g} a_j x^j$, it follows that 
$ \Delta_V(t)= (P(t))^2$, where $$P(t) = \sum_{j=0}^{g} (1 -t)^{2g - 2j} a_j
t^j.$$ From this description of $P$ if follows that: 1) $P(t)$ is symmetric, 2)
for $k \le g$, the coefficient of
$t^k$ in
$P(t)$ is a linear function of $\{a_j\}_{j = 0 \ldots k}$, and 3) in that linear
function of  
$\{a_j\}_{j = 0 \ldots k}$,
$a_k$ appears with coefficient 1.  

It follows from these observations that we can choose the $a_j$  so that $P(t) =
\Delta(t)$.  (Solve first to find $a_0 = d_0$ and then solve recursively for the
remaining $a_j$ in order.)  We now have
$\Delta_V(t) =
\Delta(t)^2$. Since    $\Delta(1) = 1$, it follows that: 4)
$a_g = 1$. Also, since $d_0 \ne 0$ we have: 5) $a_0 \ne 0$. We must now find a
matrix
$A_g$ having the desired
$\lambda$; the matrix we use  is of the following form, presented here in the
case $g = 4$.

$$\left( \begin{tabular}{c c c c}
 0 & 0 & 0 & $a_0$ \\ $-1$ & 0 & 0 & $a_1$ \\ 0 & $-1$ & 0  & $a_2$ \\ 0 & 0 &
$-1$ &
$a_3$ \\
\end{tabular}  \right)
$$

We now consider the  matrix $(1 - \omega) V + (1-\bar{\omega})V^T$.  Performing
appropriate row operations on the top $2g$ rows, and simultaneous conjugate
column operations of the first $2g$ columns, quickly yields the following matrix,
where
$\Omega = (1-\omega)(1 - \bar{\omega}) = (1-\omega) +(1 - \bar{\omega})$. (Notice
that $\Omega$ is nonzero; it equals 0 only if $\omega = 1$, which is outside our
domain.)

$$ 
\left(
\begin{tabular}{ c   c}
 {\huge{0}}$_{2g}$& 

$\left(  \begin{tabular}{c   c}                    0$_g$ & $\Omega A_g - I_g$ \\ 
 I$_g$ &  $\frac{1}{1 - \omega}I_g$ \\ \end{tabular}
\right) $

 \\  

$\left(  \begin{tabular}{c   c}                   
$0_g$ & I$_g$ \\ 
 $\Omega A_g^T -  I_g $ & $\frac{1}{1 - \bar{\omega}}I_g$ \\ \end{tabular}
\right) $

 &   \mbox{\LARGE $  \Omega    B_{2g}$} \\
\end{tabular}
\right)
$$

Choose $B_{2g}$ so that all the entries that are not in the lower right $g \times
g$ block, denoted
$B_g$, are zero.  The $\frac{1}{1  - \omega} I_g$ block can be cleared using
column operations, and simultaneous row operations will clear the $\frac{1}{1  -
\bar{\omega}}I_g$ block. It then follows that the signature of
$V_\omega$ is the signature of the following matrix.
$$\left(  \begin{tabular}{c   c}                    0 & $\Omega A_g - I_g$   \\ 
  $\Omega A_g^T - I_g$ & $\Omega$B$_g$ \\ \end{tabular}
\right) $$

Next, simultaneous column operations on the last $g$ columns and row operations
on the last $g $ rows can be used to put the $\Omega A - I_g$ block into lower
triangular form and the $\Omega A^T - I_g$ block into upper triangular form.  If
this is done, all entries on the diagonal of the upper right hand $g \times g$
block are nonzero (actually
$-1$) except the last diagonal element, which becomes
$\Omega^{g}\lambda(-\frac{1}{\Omega})$.  If $B_g$ is chosen so that all entries
are 0 except its top right and bottom left entries, let's call them $b_1$, and
the bottom right entry, say $b_2$, then after these row and column operations,
the bottom right entry of the entire matrix has become
$\Omega_2 b_2+ 2\Omega^2 a_0 b_1$.  Hence, the signature of the original matrix
$V_\omega$ is equal to the signature of the $2
\times 2$ matrix,

$$\left(  \begin{tabular}{c   c}                    0 &
$\lambda(-\frac{1}{\Omega})$   \\ 
   $\lambda(-\frac{1}{\Omega})$  & $\Omega  b_2+ 2\Omega^2 a_0 b_1$ \\
\end{tabular}
\right) $$

From the identity $ \Delta (t) =  (1-t)^{2g} \lambda(\frac{t}{(1-t)^2} )$ we see
that
$\lambda(-\frac{1}{\Omega}) = (1- \omega)^{2g}\Delta(\omega)$.   Hence  the
matrix is nonsingular with 0 signature unless $\omega$ is  root of the $\Delta$. 
On the other hand, if $\omega$ is a root of
$\Delta$, then the signature is given by the sign of   $\Omega  b_2+ 2\Omega^2
a_0 b_1$, which is the same as the sign of  $   b_2+ 2\Omega  a_0 b_1$ since
$\Omega = (1 - \omega)(1 - \bar{\omega})$ is a nonzero norm.

It remains to select $b_1$ and $b_2$ appropriately.  Suppose that $\Delta$ has
unit roots with positive imaginary parts $\{\omega_1, \ldots , \omega_k\}$,
placed in order of increasing real part, and let $p$ satisfy $1 \le p \le k$. 
Since $\Omega = 2 - 2\mbox{Re\;}\omega$  is a decreasing function of the real
part of
$\omega$, it is clear that $b_1$ and $b_2$ can be selected so that $   b_2+
2\Omega  a_0 b_1$ is positive for $\omega_j$ when $ j \le p$ and it is negative
for $\omega_j$ when $ j >  p$. The corresponding matrix, say $V_1$, has signature
$\sigma_{\omega_j}(V_1) = 1$ if $k \le p$ and
$\sigma_{\omega_j}(V_1) = -1$ if $k > p$.  Similarly, construct $V_2$ so that
$\sigma_{\omega_j}(V_2) = -1$ if $k < p$ and
$\sigma_{\omega_j}(V_2) =  1$ if $k \ge p$. The direct sum of these two matrices
has the desired property: its only nontrivial signature is  t
$\sigma_{\omega_p}(V_1 \oplus V_2) = 2$.

\end{proof}

\noindent{\bf Comment.}  For the matrix constructed we have $\sigma_{\omega_p} =
2$ and all other signatures are 0.  Can a similar matrix be constructed, only
with $\sigma_{\omega_p} = 1$?  If the irreducible polynomial for $\omega_p$ has a
unique root $\omega \in S$ the construction gives such an example.  If, however,
there are multiple roots in $S$, it can be shown that the parity of the
signatures at each of these roots must be the same, as follows.  The parity of
the signature of
$V_\omega$ is determined by its rank and nullity of $V_\omega$.  Both of these
are invariant under the action of the Galois group permuting the roots of
$\Delta(t)$.
 %%%%%%%%%%%%%%%%SECTION%%%%%%%%%%%

\section{Higher Dimensional Knots} The algebraic theory of 1 dimensional knots in
$S^3$ extends to a general theory of codimension 2 knots in $S^{2n +1}$. 
According to Levine \cite{le1}, the case of knots in dimension $S^{4n + 3}$ is
identical to the classical case.  Hence, all  the results presented so far apply
for knots in these dimensions.  In this section we will describe how to modify
our previous work to apply to knots in $S^{4n+1}$. The reference for this is
\cite{le1}.

  In the case of knots in $S^{4n +1}$, a Seifert matrix is a $2g \times 2g$
integral matrix satisfying $\det( V + V^T) = 1$.  The signature function of such
a Seifert matrix is given by
$$
\sigma_\omega(V)=\mathrm{sign}\big[(\omega-\bar\omega)\left((1-\omega)
V-(1-\bar\omega)V^T\right)\big].
$$  The Alexander polynomial is given by $\Delta(t)=\det(tV+V^T)$.  An integral
polynomial is the Alexander polynomial of some such Seifert matrix if and only if
$\Delta(t)=t^{2g}\Delta(t^{-1})$, $\Delta(1)=(-1)^g$, and
$\Delta(-1)$ is a square.

\vskip.1in
\noindent{\bf Extending Theorem~\ref{thmind}.}  To extend the proof of
Theorem \ref{thmind} to the case of knots in $S^{4n+1}$, we need to modify the
  polynomial $\Delta(t)$ constructed in the proof to assure that it
satisfies the stricter conditions on Alexander polynomials in these dimensions. 
To do this, we replace $\Delta(t)$ with the polynomial 
$D(t)=(ct^2+(1-2c)t+c) \Delta(t)$ where $c=2(a+b)>0$.  Then it is straightforward
to check that $D(t)$ is an Alexander polynomial of a Seifert matrix $V$ of a
knot in $S^{4n+1}$, using the fact that
$\Delta(1)=-1$ and $\Delta(-1)=8a+8b-1$.  Even though $D(t)$ has an additional
zero $\omega''$ in $S$, where
$\mathrm{Re}(\omega'')=1-1/2c$, we can control it by choosing $a$ and
$b$ carefully, as follows.  Since $a/b$ is to be  chosen close to $r > -1$, we can
assume that $a/b > -1 + \epsilon$ for some fixed positive $\epsilon$.  Hence, 
$c=2(a+b)>2b\epsilon$ and  we can choose 
$b$ large enough so that
$1-1/2c$ is sufficiently close~$1$.  Since both $\sigma_1$ and
$\sigma_{(-1)}$ are   zero for any knot in $S^{4k+1}$, it follows
that the signature function $\sigma_\omega(V)$ assumes a constant nonzero value
(indeed $\pm 2$) for
$r<\mathrm{Re}(\omega)<1-1/2c$, and zero for $\mathrm{Re}(\omega)<r$ or
$\mathrm{Re}(\omega)>1-1/2c$.  Thus it can be used to complete the proof in this
setting.

\vskip.1in
\noindent{\bf Extending Theorem~\ref{thmmet}.} The extension of the proof of
Theorem~\ref{thmmet} for knots in 
$S^{4n+1}$ is more straightforward.  We proceed in the same way
but start with
$$
V=\left(\begin{matrix}
\mbox{\LARGE $0_{2g}$} &
\left(\begin{matrix} I_g & A_g \\ I_g & 0_g \end{matrix}\right) \\
\left(\begin{matrix} 0_g & -I_g \\ -A_g^T & I_g \end{matrix}\right) &
\mbox{\LARGE $B_{2g}$}
\end{matrix}\right),
$$ 
which is a Seifert matrix of a   slice knot in $S^{4n +1}$ for
any $n$ by a result of~\cite{le1}, and assume that the given
$\Delta(t)$ satisfies $\Delta(1)=(-1)^g$.  Then it can be shown that
$\sigma_\omega(V)$ is zero if $\Delta(\omega)\ne 0$, and
$\sigma_\omega(V)=b_2-2\Omega a_0 b_1$ if $\Delta(\omega)= 0$.  Now
the arguments of the last paragraph of the proof of Theorem~\ref{thmmet} can be
applied to construct a desired matrix.

\section{Applications to Knot Concordance}\label{appsect}

Our detailed examination of signature functions was motivated by problems in studying
the concordance group of knots.  In this section we give some indication as to the
usefulness of the algebraic results of this paper.

The Seifert form $V_K$ associated to a knot $K$ in $S^3$ is not unique.  However, as
mentioned in Section~\ref{concord}, by placing an equivalence relation on the set of
knots one arrives at the {\it concordance group} of knots, $\calc$.  The association
$K \to V_K$ induces a surjective homomorphism: $\phi: \calc \to \calg$.  These
notions were defined and studied by Levine \cite{le1, le2}.

Levine's work  contained two main results.  One was topological---in
higher dimensions the analog of $\phi$ is an isomorphism.  The second was
algebraic---the group $\calg$ is isomorphic to an infinite direct sum, $\zz^\infty
\oplus \zz_2^\infty \oplus \zz_4^\infty$.
This algebraic result depended in part on the existence of an infinite
collection of $\omega$ for which the associated signature functions $\sigma_\omega$
are linearly independent.

Several years after Levine's work, Casson and Gordon \cite{cg1, cg2} proved that
$\phi$ is not an isomorphism (for knots in $S^3$) by developing obstructions to a
knot with metabolic Seifert form from being
trivial in
$\calc$. The kernel of $\phi$ is called the group of {\it algbebraically slice}
knots, denoted $\cala$.  Later, Gilmer
\cite{g1, g2} demonstrated that these Casson-Gordon obstructions could be
interpreted in terms of knot signatures: not those of the original knot, which are
necessarily 0 if the knot is algebraically slice, but rather a knot that
reflects the knotting in a surface bounded by the original knot.

The most basic applications of Gilmer's approach to Casson-Gordon invariants called
on rather simple facts about the signature function.  For instance, constructing
algebraically slice knots that are nontrivial in the concordance group was reduced
to finding knots with nontrivial signature at a single root of unity.

As the subtleties of $\cala$ have been explored, deeper facts about signatures have
been called on. As one example, Stoltzfus \cite{st}  proved that in higher
dimensions if  the Alexander polynomial of a knot $K$ factors into irreducible
factors with resultant 1, then the associated knot is concordant to a corresponding
connected sum of knots.  The second author of this paper, in unpublished work, has
shown that this result does not apply in dimension three.  The simplest example
involves the polynomials $2t^2 - 3t +2$ and $3t^2 - 5t +3$.  The construction of
the counterexample depended, via Gilmer's work,  on finding a knot whose signatures
at 67--roots of unity satisfy a complicated linear inequality.  Finding such a knot
could be carried out by ad hoc methods, but  Theorem~\ref{thmind} makes the
existence of such a knot automatic.

  If one takes on more general problems, the ad hoc methods that can be applied to a
single example are no longer useful.  For instance, it now appears to be the
case that for almost any pair of Alexander polynomials with resultant 1 there is a
knot with Alexander polynomial the product of those polynomials, and yet the knot is
not concordant to a corresponding connected sum.  The proof depends on contructing
knots whose signature function at the collection of some (unknown) roots of unity
satisfy some (unknown) inequality.  Because of the general nature of the problem,
little about which roots of unity or what the inequality is  can be identified. 
Yet, with Theorem~\ref{thmind} it is possible to assert that such a knot will exist.

Similar issues arise in a number of related settings.  In finding properties of
$\calg$ that do not apply to $\calc$, individual examples can sometimes be
constructed using (perhaps very messy) ad hoc constructions, but general results
demand complete control over the signature function, something that is offered by
Theorem~\ref{thmind}.  

Applications of Theorem~\ref{thmmet} take place in a different realm.  Levine's work
in \cite{le3} offered one proof that the quotient $\calg / \calg^h$ is nontrivial,
where
$\calg^h$ is the algebraic concordance group built using the equivalence relation
based on hyperbolic rather than metabolic forms. In fact, it follows from \cite{le3}
that the quotient is infinitely generated.  This result implied a  similar result for
the topological setting of double null concordance versus concordance. Levine's work
focused on the signature function at a particular root of unity.   Recent work of Cochran, Orr, and Teichner~\cite{cot} has
renewed interest in the study of double null concordance; in particular, Taehee Kim,
in unpublished work, has made significant progress in studying the case in which all
invariants of the knot based on $\calg$ and $\calg^h$ vanish.  This ongoing work
points to a need for a further study of the algebra of $\calg^h$.  The work here
demonstrates the a complete analysis of
$\calg / \calg^h$ will depend on considerations of all possible unit roots of
Alexander polynomials.  
%%%%%%%%%%%%%REFERENCES

\newcommand{\etalchar}[1]{$^{#1}$}

\end{document}